\documentclass{article}
\usepackage{amsthm, amsmath, amssymb, amsfonts}

\def\RR{\mathbb{R}}
\def\Trop{\mathrm{Trop\ }}

\newtheorem*{Theorem}{Theorem}
\newtheorem*{Lemma}{Lemma}

\theoremstyle{definition}

\newtheorem{problem}{Problem}

\begin{document}

\title{Reconstructing Trees from Subtree Weights}
\author{Lior Pachter and David Speyer\\ Department of Mathematics, U.C. Berkeley, Berkeley CA 94720}
\maketitle
\begin{abstract}
The tree-metric theorem provides a necessary and sufficient condition for a dissimilarity matrix to be a tree metric, and has served as the foundation for numerous distance-based reconstruction methods in phylogenetics. Our main result is an extension of the tree-metric theorem to
more general dissimilarity maps. In particular, we show that a tree with $n$ leaves is reconstructible from the weights of the $m$-leaf subtrees provided that $n \geq 2m-1$. 
\end{abstract}

\section{Introduction}

The problem of reconstructing a graph from measurements of distances between certain nodes was first proposed by Hakimi and Yau \cite{Hakimi} in 1965, and has since developed into a myriad of related questions arising from diverse applications such as phylogeny reconstruction and internet tomography. These problems all have in common the notion of a graph realization of a matrix. That is, if $D$ is a matrix whose rows and columns are indexed by a set $X$, then 
$D$ has a {\it realization} if there is a weighted graph $G$ whose node set contains $X$, and that satisfies $d(u,v)=D(u,v)$ where $d(u,v)$ is the distance between two nodes in the graph. In what follows we will assume that our matrices $D$ are symmetric and have zeros on the diagonal; in phylogenetics these are called {\it dissimilarity matrices} \cite{Semple}. A non-negative dissimilarity matrix that has a graph realization is called a {\it distance matrix}. The Hakimi-Yau problem is therefore:
\begin{problem}[Hakimi-Yau]
Given a distance matrix on a set $X$, find a graph $G$ which is a realization of $D$ so that the sum of the edge lengths of $G$ is minimized.
\end{problem}
The general problem is notoriously difficult \cite{Chung}, however the case when $G$ is restricted to be a tree (and $X$ corresponds to leaves of the tree) has been well understood. Distance matrices realized by trees are said to be {\it tree metrics}, and the main result is the classic {\it tree-metric theorem} \cite{Buneman1, Buneman2,Simoes,Zaretskii}: 
\begin{Theorem}[Tree-Metric Theorem]
\label{thm:tree-metric}
Let $D$ be non-negative dissimilarity matrix on $X$ . Then $D$ is a tree metric on $X$ if and only if, for every four (not necessarily distinct) elements $i,j,k,l \in X$, two of the three terms in the following equation are equal and greater than or equal to the third:
\begin{equation}
\label{def:four-point}
D(i,j)+D(k,l),\ D(i,k)+D(j,l),\ D(i,l)+D(j,k) 
\end{equation}
Furthermore, the tree $T$ with leaves $X$ that realizes $D$ is unique. 
\end{Theorem}
The {\it four point condition} (\ref{def:four-point}) is therefore a necessary and sufficient condition for 
a matrix to be realized by a tree. 

The algorithmic question of how to construct a tree from a tree metric is particularly important for the reconstruction of phylogenetic trees from genomic sequences, because the tree-metric theorem suggests the approach of estimating pairwise distances between the sequences and then reconstructing the tree. Methods that reconstruct the correct unique tree given a tree-metric, are called {\it consistent}, and biologists are particularly interested in consistent methods that are able to reconstruct the correct tree topology even with inaccurate distance matrices. Numerous polynomial time consistent algorithms have been proposed, of which we mention the classic Buneman method \cite{Buneman1} and the neighbor-joining method \cite{Saitou}. The neighbor-joining algorithm (NJ) is very popular because it is fast and performs well in practice; we mention the Buneman construction because the proof of our main result is based on it.

In this paper we propose the following generalization of the Hakimi-Yau problem: An {\it $m$-dissimilarity map} is a map $D:X^m \rightarrow {\mathbb R}$ ($X^m$ denotes the $m$th Cartesian product of $X$) with $D(x_1,\ldots,x_m)=D(x_{\pi(1)},\ldots,x_{\pi(m)})$ for all permutations $\pi \in S_m$ and $D(x,x,\ldots,x)=0$. 
We say that a graph $G$ realizes $D$ if the node set of $G$ contains $X$ and for every $x_1,\ldots,x_m\in X$, the weight of the smallest subgraph in $G$ containing $x_1, \ldots, x_m$ is $D(x_1,\ldots,x_m)$. We call an $m$-dissimilarity map that is realizable an {\it $m$-distance map}.

\begin{problem}
Given an m-distance map $D$ on a set $X$, find a graph $G$ which is a realization of $D$ so that the sum of the edge lengths of $G$ is minimized.
\end{problem}
 
Our main result is a weak version of the tree-metric theorem for $m$-dissimilarity maps: we show that an $m$-dissimilarity map does determine the tree it comes from (as long as $m$ is small compared to the number of leaves) although we give no explicit criterion for an $m$-dissimilarity map to come from a tree. We also discuss the ramifications of our result for phylogenetic tree reconstruction, in particular the relationship between our reconstruction theorem and maximum-likelihood methods for phylogenetic trees.

\section{Reconstructing Trees}

Let $T$ be a (finite) tree, with a positive weight $w(e)$ assigned to each edge $e$. Let $L(T)$ denote the set of leaves of $T$.  For any $V \subseteq L(T)$ of leaves, let $[V]$ denote the smallest subtree of $T$ containing $V$. For $T'$ any subtree of $T$, let $w(T')$ be the sum of the weights of the edges of $T'$.  
The main result of this paper is 
\begin{Theorem}
Let $T$ be a tree with $n$ leaves and no vertices of degree 2. Let $m \geq 3$ be an integer. If $n \geq 2m-1$, then $T$ is determined by the set of values $w([V])$ as $V$ ranges over all $m$ element subsets of $L(T)$. If $n=2m-2$, this is not true.
\end{Theorem}
In other words, if $m$ is fixed and the number of leaves in a tree $T$ is sufficiently large, then 
$T$ is uniquely determined from the $m$-dissimilarity map corresponding to the weights of the subtrees on $m$ leaves. 

We will engage in the following abuse of notation: when $S \subset L(T)$ and $v_1$,\dots, $v_r \in L(T)$, we will write $v_1 \ldots v_r$ to denote $\{v_1,\ldots, v_r\}$ and we will write $S v_1\ldots v_r$ to mean $S \cup \{ v_1, \ldots, v_r \}$. We will also write $w(V)$ to mean $w([V])$.

Let $T$ be a tree and let $i$, $j$, $k$ and $l$ be distinct members of $L(T)$. 
We will say $i$ and $j$ \emph{are opposite} $k$ and $l$ if $[ijkl]$ is of the following form: there are two (distinct) vertices $v$ and $w$ of $T$ such that $[ijkl]$ consists of five edge distinct paths: a path from $i$ to $v$, a path from $j$ to $v$, a path from $v$ to $w$, a path from $w$ to $k$ and a path from $w$ to $l$. We will denote this state of affairs by $(i,j;k,l)$ and we will denote the path from $v$ to $w$ by $\gamma_{ijkl}$. It is clear that $\gamma_{ijkl}$ is well defined; in phylogenetics it is known as the Buneman index of the split $(i,j;k,l)$. 
Finally, if $S$ is a set, $\binom{S}{m}$ denotes the set of $m$ element subsets of $S$.

In order to prove the main theorem, we first show that if $T$ is a tree with no vertices of degree 2 and having $n$ leaves ($n \geq 2m-1$, $m \geq 3$) then we can reconstruct the topology of $T$. We will rely on the well known result that $T$ can be recovered as an abstract graph (without determining the edge weights) by knowing the signs of the Buneman indices, i.e. for which $i$, $j$, $k$ and $l$ we have $(i,j;k,l)$. The following Lemma that shows that this data can be recovered from the values $w([V])$ as $V$ ranges over $\binom{L(T)}{m}$. 

\begin{Lemma}
Let $i$, $j$, $k$ and $l$ be distinct members of $L(T)$. Then $(i,j;k,l)$ iff there is an $R \in \binom{L(T) \setminus \{i,j,k,l\}}{m-2}$ such that 
\begin{equation}
\label{eq:mtree}
w(Rij)+w(Rkl) < w(Rik)+w(Rjl)=w(Ril)+w(Rjk). \qquad 
\end{equation}
\end{Lemma}

\begin{proof}
Let $R \in \binom{L(T) \setminus \{i,j,k,l\}}{m-2}$. Let $\tilde{T}$ be the tree that arises by contracting $[R]$ to a point. We will denote the images of $i$, $j$, etc. in $T$ by $\tilde{i}$, $\tilde{j}$, etc. We note that $\tilde{R}$ is a leaf of $\tilde{T}$. Now,
$$w(Rij)=w(\tilde{R}\tilde{i}\tilde{j})+w(R),$$
so (\ref{eq:mtree}) holds iff the corresponding statement holds in $\tilde{T}$.

If $(\tilde{i}, \tilde{j}; \tilde{k}, \tilde{l})$ then clearly $(i,j;k,l)$. We claim that, conversely, if $(i,j;k,l)$ then there exists some $R$ such that $(\tilde{i}, \tilde{j}; \tilde{k}, \tilde{l})$. To see this, let $e$ be an edge in $\gamma_{ijkl}$. Removing $e$ from $T$ divides $T$ into two trees, let $L_1$ and $L_2$ be the leaves of these trees. Now, $|L(T) \setminus \{i,j,k,l\}|=n-4 \geq 2m-5$ so either $L_1 \setminus \{i,j,k,l\}$ or $L_2 \setminus \{i,j,k,l\}$ has at least $m-2$ members, without loss of generality suppose the former does. Then take $R \in \binom{L_1 \setminus \{i,j,k,l\}}{m-2}$ and you will have $e \not\in R$ so $(\tilde{i}, \tilde{j}; \tilde{k}, \tilde{l})$.

The previous two paragraphs reduce us to the case where $|R|=1$, that is, where $m=3$ and $R=\{ r \}$ for some $r \in L(T)$. We may delete all leaves other than $i$, $j$, $k$, $l$ and $r$. Now the claim is reduced to the case where $m=3$ and $n=5$, and can be verified by checking it for all trees with five leaves.
\end{proof}

\textbf{Remark:} 
One can check that, if the subtree weights are modified by less that $e/2$, where $e$ is the weight of the smallest edge, the reconstructed tree will not be altered. This bound of $e/2$ is the same as for the classical case of $m=2$ \cite{Atteson}.

We still need to show that we can recover the edge weights of $T$.
Continue to assume $T$ has no vertices of degree 2, $|L(T)|=n$ and $n \geq 2m-1$. By the Lemma, the values $w(V)$, $V \in \binom{L(T)}{m}$, determine the graph of $T$, so we may assume that this is known. 

Let $e$ be an edge of $T$. We first consider the case in which neither of the endpoints of $e$ is a leaf. Then we can find $i$, $j$, $k$ and $l$ in $L(T)$ such that $(i,j;k,l)$ and $\gamma_{ijkl}=e$ (here and at many points in the future, we use that we already know the graph of $T$). Moreover, using $n-4 \geq 2m-5$ as before, we can find $R \in \binom{L(T) \setminus \{i,j,k,l\}}{m-2}$ such that $ e \not\in [R]$. Then
$$w(Rik)+w(Rjl)-w(Rij)-w(Rkl)=2w(e)$$
so we can determine $w(e)$.

Let $v \in L(T)$, let $e_v$ be the unique edge incident to $v$. Let us assume that we have already computed $w(e)$ for all $e \in T$ not of the form $e_v$. Then, for all $V \subset \binom{L(T)}{m}$, we can compute 
$$\sum_{v \in V} w(e_v)=w(V)-\sum_{\substack{e \in V \\ e \neq e_v}} w(e).$$

Now, let $i$ and $j$ be distinct members of $L(T)$. For any $S \in \binom{L(T) \setminus \{ i, j \}}{m-1}$  we have 
$$\sum_{v \in Si} w(e_v)-\sum_{v \in Sj} w(e_v)=w(e_i)-w(e_j).$$
Thus, we can determine the $w(e_v)$ up to a common additive constant. We may determine that constant from $w(V)$ for any $V \in \binom{L(T)}{m}$.

It remains to show that if $n=2m-2$ it may not be possible to reconstruct the tree.
Let $T$ be the following tree: the vertices of $T$ are known as $v_1$, \dots, $v_n$ and $w_2$, \dots $w_{n-1}$. The edges of $T$ are of the following two forms: $(v_{i-1}, v_{i})$ and $(v_i, w_i)$. We will set $w_1=v_1$ and $w_n=v_n$. Clearly, $L(T)=\{ w_1=v_1, w_2, \dots, w_{n-1}, w_n=v_n \}$. The edge weights may be chosen arbitrarily.

Let $T'$ be the same tree except that the edges $(v_{m-1}, w_{m-1})$ and $(v_m, w_m)$ are deleted and replaced by $(v_{m}, w_{m-1})$ and $(v_{m-1}, w_{m})$. We place the edges of $T$ and those of $T'$ in bijection by making $v_{m-1} w_{m-1}$ correspond to $v_m w_{m-1}$, $v_m w_m$ correspond to $v_{m-1} w_m$ and pairing all other edges in the obvious way. Assign the weights to the edges of $T'$ that are assigned to the corresponding edges of $T$.

We claim that, for any $V \in \binom{L(T)}{m}=\binom{L(T')}{m}$, the same edges (using the above bijection) appear in $[V]$ and $[V]'$, where $[V]'$ denotes the minimal subtree of $T'$ containing $V$. To prove this, we consider divide the edges $e$ of $T$ into three types. 

First, we could have $e=(v_i,w_i)$. Then $e \in [V]$ iff $w_i \in V$ and the same is true for $[V]'$. Secondly, we could have $e=(v_{i-1}, v_{i})$ with $i \neq m$. Then $e \in [V]$ iff $w_j$ and $w_k \in V$ for some $j \leq i < i+1 \leq k$. This also holds for $[V]'$.

So far, we have used the fact that $|V|=m$, the previous paragraph is correct for any $V \subseteq L(T)$. We now consider the final case, $e=(v_{m-1}, v_m)$. Removing $e$ from $T$ divides $L(T)$ into two sets $L_1$ and $L_2$, each of size $m-1$. As $|V|=m$, $V \cap L_1$ and $V \cap L_2 \neq \emptyset$ so $e \in [V]$. Similarly, $e \in [V]'$.
\qed

Our main result shows that, when $n$ is large enough compared to $m$, we can reconstruct a tree from the weights of its $m$-leaf subtrees. However, if we are simply given an $m$-dissimilarity map $D: \binom{[n]}{m} \to \RR$, we do not know how to test whether this map comes from a tree.

When $m=2$, this is given by the tree metric theorem. When $m$ is larger than $2$, an obviously necessary condition is that, for every $R \in \binom{[n]}{m-2}$ and $i$, $j$, $k$ and $l \in [n] \setminus R$, two out of three of the following expressions must be equal to each other and greater than or equal to the third:
$$D(Rij)+D(Rkl),\ D(Rik)+D(Rjl),\ D(Ril)+D(Rjk) .$$
Moreover, we can impose the combinatorial requirement that when we consider the above equation with the same $(i,j,k,l)$ and different $R$ and $R'$, that the same one of the three terms above is minimized. However, by counting dimensions, we can see that this condition is not adequate in any case except $n=5$, $m=3$.

\section{The Tropical Analogy}

In this section, we describe a connection between subtree weights and an area of algebraic geometry known as ``tropical geometry'', and inquire whether the analogy can be made tighter. The basic reference for our discussions is \cite{SS}, although we reverse the sign conventions of that paper to more closely match those occurring elsewhere in this one.

Let $f=\sum_{e \in E} f_{e_1 \cdots e_n} x_1^{e_1} \cdots x_n^{e_n}$ be a polynomial in $n$ variables, where the $f_{e_1 \cdots e_n}$ are nonzero. Define $\Trop f$ to be the subset of $w=(w_1, \ldots, w_n) \in \RR^n$ such that, of the collection of numbers $\sum_{i=1}^n e_i w_i$ where $e$ runs over $E$, the maximum occurs twice. If $I \in K[x_1, \ldots, x_n]$ is an ideal, set $\Trop I=\bigcap_{f \in I} \Trop f$. 

Now, consider the case of a polynomial ring whose variables are indexed by the $m$ element subsets of $[n]$, we will write these variables as $p_{S}$ for $S \in \binom{[n]}{m}$. Let $w \in \RR^{\binom{n}{m}}$. The statement that the maximum of 
$$w_{Rij}+w_{Rkl},\ w_{Rik}+w_{Rjl},\ w_{Ril}+w_{Rjk}$$
occurs twice precisely says that
$$w \in \Trop ( p_{Rij} p_{Rkl} - p_{Rik} p_{Rjl} + p_{Ril} p_{Rjk}).$$

So, if $w$ arises from the $m$-leaf subtree weights of a tree, the $w \in \Trop f$, where $f$ is any polynomial of the form $p_{Rij} p_{Rkl} - p_{Rik} p_{Rjl} + p_{Ril} p_{Rjk}$, $i < j < k < l$. Such polynomials are called the \emph{three term Plucker relations}. It is well known that all of these relations lie in the ideal of the Grassmanian $G(m,n)$. 

In the case where $m=2$, it was shown in \cite{SS} that $\Trop G(2,n)$ is exactly the space of tree metrics. This raises several natural problems:

\begin{problem}
Does the space of $m$-leaf subtree weights lie in $\Trop G(m,n)$?
\end{problem}

(It can be shown that the two spaces are not equal: $\Trop G(m,n)$ has larger dimension.) Assuming the answer to the above problem is ``yes'', and inspired by the bijection between tree metrics and points of $\Trop G(2,n)$, we can ask for more.

\begin{problem}
Is there a map $G(2,n) \to G(m,n)$ with image $X$ such that $\Trop X$ is the space of $m$-leaf subtree weights?
\end{problem}

A positive answer to the above problem is known in the case where $m=3$. Write a point of $G(m,n)$ as a matrix with $n$ rows and $m$-columns, considered up to the right action of $GL_n$. One can check that the following map $Mat_{2 \times n} \to Mat_{3 \times n}$ descends to a map $G(2,n) \to G(3,n)$:
$$\begin{pmatrix}
x_{11} & x_{12} & \cdots & x_{1n} \\
x_{21} & x_{22} & \cdots & x_{2n} \end{pmatrix} \mapsto
\begin{pmatrix}
x_{11}^2 & x_{12}^2 & \cdots & x_{1n}^2 \\
x_{11} x_{21} & x_{12} x_{22} & \cdots & x_{1n} x_{2n} \\
x_{21} & x_{22} & \cdots & x_{2n} \end{pmatrix}$$

This map takes a tree metric to twice the corresponding $3$-leaf subtree weight. We have not found such a ``geometric lifting'' for $m>3$.

\section{Applications}
The fundamental problem in phylogenetics is to reconstruct trees from sets of sequences related by an evolutionary tree. The sequences can be DNA or protein sequences, or more generally can encode the order of genes in a genome or other evolving features. Some of the most popular methods for reconstructing trees are {\it distance based methods}. Distance based methods 
begin by estimating pairwise distances between the sequences, thus leading to a dissimilarity map (although it need not be an actual tree-metric). A tree is then reconstructed from the dissimilarity map, and it is hoped that the topology of the tree is correct.

Distance-based method are typically based on maximum-likelihood estimates of the pairwise distances. Suppose that the sequences are labeled $s^1,\ldots,s^n$ and that position $k$ in the $j$th sequence is denoted by $s_k^j$. For a pair of sequences $i,j$ we assume that one of them, $s^j$ has a the $k$ position picked from a distribution $q_{s_k}^i$, and that the other is related by
a multiplicative substitution process to the first sequence, i.e. the $k$th character in 
$s^i$ is distributed according to $P(s_k^i|s_k^j,t)$ where $\sum_b P(s_k^i|s_k^l,t_1)P(s_k^l|s_k^j,t_2)=P(s_kk^i|s_k^j,t_1+t_2)$. 
Furthermore, we assume that the substitution process is reversible so that we can switch the indices $i,j$. 

Then the pairwise distances are computed by
\[ d_{ij} = argmax_{t} \left(  \prod_k q_{s_k}^j P(s_k^i|s_k^j,t) \right). \]
Since the terms $q_{s_k}^j$ do not depend on $t$ we can write this as
\[ d_{ij} = argmax_{t} \left( \prod_k P(s_k^i|s_k^j,t) \right). \]
When the number of sites is large the consistency of maximum likelihood implies that these estimates will converge to the true branch lengths in the tree (assuming that the probabilistic model is correct). Thus, distance methods based on ML pairwise distance estimates can be viewed as generating approximations of the actual ML tree.  

In practice, the use of short sequences can lead to inaccurate estimates, especially for longer branch lengths, and this is a common source of error for methods such as neighbor joining \cite{Saitou}. A number of solutions to this problem have been proposed: for example variants of neighbor joining exist that are less sensitive to long branch lengths \cite{Bruno}. Ranwez and Gascuel have shown that long pairwise branch lengths can be corrected  by considering ML estimates of distances using three taxa \cite{Ranwez}. Quartet methods such as quartet puzzling \cite{Strimmer} attempt to reconstruct the tree from more reliable quartets rather than pairwise distances. However a major problem with quartet methods has been how to reconcile the different topologies of the quartets into one tree.

Our approach can be seen as generating a better approximation to the ML tree, by relying on more accurate $m$-tree weights rather than pairwise distance. By avoiding the need to reconcile diverse topologies, we avoid the difficulties of standard quartet methods. Furthermore, as $m$ increases, the splits can be more accurately identified and thus a more accurate tree reconstructed. Although the computational complexity of subtree weight reconstruction grows rapidly with $m$, it remains polynomial, and furthermore many key steps are trivially parallelizable. 

\section*{Acknowledgments}

We thank Bernd Sturmfels for many comments which improved the manuscript. Lior Pachter was partially supported by
a grant from the NIH (R01-HG02362-02).

\end{document}